\documentstyle[12pt]{article}
\newcommand{\vs}{\vspace{.1in}}

\newcommand{\noin}{\noindent}

\begin{document}

\begin{center} {\bf Non-Hyperbolic Complex Spaces with Hyperbolic Normalizations} \end{center} 

\begin{center} Shulim Kaliman \\

Department of Mathematics and Computer Science \\

University of Miami \\

Coral Gables, FL  33124, U.S.A. \end{center}

\begin{center} Mikhail Zaidenberg \\

Universit$\acute{e}$ Grenoble I. \\

Istitut Fairier de Math$\acute{e}$matiques \\

38402 St. Martin d'H$\grave{e}$res-Cedex \\

France \end{center}

\begin{center} {\bf 1. Introduction.} \end{center}

Let $E$ be a Hausdorff topological space and $O$ be a subsheaf of the sheaf 
of germs of continuous functions on $E$. Recall that the pair $(E, O)$ is 
called a complex space, if every point of $E$ has an open neighborhood $U$ 
such that one can identify $U$ with an analytic subset of an open subset 
$\Omega \subset {\bf C}^{n}$ and $O\mid _{U}$ with the 
sheaf of germs of holomorphic functions on the analytic set $U$. As usual we shall denote 
$(E, O)$ by $E$. A complex space $E$ is said to be normal at a point 
$a\in E$, if the ring of germs of holomorphic functions at $a$ is 
integrally closed in its quotient field. If $E$ is normal at its every point, then 
$E$ is said to be normal. As it is well-known every smooth complex manifold 
is normal. A normalization of a complex space $E$ is a pair $(\tilde{E} , 
\pi )$ consisting of a normal complex space $\tilde{E}$ and 
a surjective holomorphic mapping $\pi : \tilde{E} \rightarrow E$ such that \vs

$(i)$ $\; \pi :\tilde{E} \rightarrow E$ is proper and $\pi ^{-1} (a)$ is 
finite for every $a\in E$; \vs

$(ii)$ if $S$ is the set of singular points of $E$, then $\tilde{E} -\pi 
^{-1} (S)$ is dense in $\tilde{E}$ and $\pi :\tilde{E} -\pi ^{-1} (S) \rightarrow E-S$ is biholomorphic. \vs

A complex space is called hyperbolic, if the Kobayashi pseudodistance on it 
is a distance. \vs

{\bf THEOREM ([2], [3]).} {\em If $(\tilde{E} , \pi )$ is a normalization 
of a hyperbolic space $E$, then $\tilde{E}$ is also hyperbolic.} \vs

One can ask if the converse statement is true. More precisely, 
if $E$ is a hyperbolic complex space, is $\tilde{E}$ hyperbolic? The answer is negative 
and we shall present two very simple examples below. \vs

The second author is grateful to Mathematisches Institut and SFB-170 
``Geometric und Analysis'' for support during the preparation of this 
paper. \vs

\begin{center} {\bf 2. Examples.} \end{center} 

Let $(x,y,u,v)$ be a coordinate system in ${\bf C}^{4}$ and $E$ be the 
affine algebraic subvariety given by the equations

\[ \left\{ \begin{array}{rl} y^{4} & =x^{4} -1 \\

u^{4} & =y^{4} (v^{4} -1). \end{array} \right. \]

\noin Then $E$ is not a hyperbolic complex space, since the set $E\cap 
\{ y=0\}$ consists of complex lines. 
Let $\tilde{E}$ be the affine algebraic submanifold in ${\bf C}^{4}$ given by the equations 

\[ \left\{ \begin{array}{rl} y^{4} & =x^{4} -1 \\

u^{4} & =v^{4} -1. \end{array} \right. \]

\noin Then $\tilde{E}$ is smooth and, in particular, normal. One can 
consider $\tilde{E}$ as the direct product $R\times R$, where $R$ is the 
smooth algebraic curve defined by the equation 
$x^{4}_{1} =x^{4}_{2} -1$ in ${\bf C}^{2} =\{ (x_{1} , x_{2} )\}$. Since $R$ is hyperbolic, $\tilde{E}$ 
is also hyperbolic ([2, Proposition 4.1]). Put $p(x,y,u,v)=(x,y,yu,v)$ and 
$\pi =p\mid _{\tilde{E}}$. Clearly, $(\tilde{E} , \pi )$ is a normalization 
of $E$. This is a desired example. \vs

Now we shall present a proper example in the compact case. Choose a 
homogeneous coordinate system $(X,Y,Z)$ in ${\bf CP}^{2}$. Put $\overline{R}
=\{ (X,Y,Z)\in {\bf CP}^{2} \mid Y^{4} =X^{4} -Z^{4} \}$. Then $\overline{R}$ 
is a smooth Riemann surface of genus 3, and, therefore, $\overline{R}$ is 
hyperbolic. Put $R_{1} =\overline{R} -\{ Z=0\}$ and 
$R_{2} =\overline{R} -\{ X=0\}$. Let $q:B\rightarrow \overline{R}$ 
be a holomorphic fiber bundle 
over $\overline{R}$ with fiber ${\bf CP}^{2}$ such that $q^{-1} (R_{k} )
\cong R_{k} \times F_{k}$, where $F_{k} \cong {\bf CP}^{2}$ with a 
homogeneous coordinate system $(U_{k} , V_{k} , W_{k} )$, and the 
connection between the coordinate systems in $q^{-1} (R_{1} \cap R_{2} )$ is given by 
the formula $(U_{2} , V_{2} , W_{2} )=(ZU_{1}, XV_{1}, XW_{1} )$. Consider 
the compact complex space $E_{1} \subset B$ so that 

\[ E_{1} \cap q^{-1} (R_{1} )= \left\{ \begin{array}{rl} Y^{4} & =X^{4} -Z^{4} \\

Z^{4} U^{4}_{1} & =Y^{4} (V^{4}_{1} -W^{4}_{1} ) \end{array} \right. \]

\noin and 

\[ E_{1} \cap q^{-1} (R_{2} )= \left\{ \begin{array}{rl} Y^{4} & =X^{4} -Z^{4} \\

X^{4} U^{4}_{2} & =Y^{4} (V^{4}_{2} -W^{4}_{2} ) \end{array} \right. \]

Then $E_{1}$ is not hyperbolic, since the set $E_{1} \cap \{ Y=0\}$ consists 
of Riemann spheres. Let $\tilde{E}_{1}$ be the complex space in $B$ given by 
the formulas

\[ \tilde{E}_{1} \cap q^{-1} (R_{1} )= \left\{ \begin{array}{rl} Y^{4} & =X^{4} -Z^{4} \\

Z^{4} U^{4}_{1} & =V^{4}_{1} -W^{4}_{1}  \end{array} \right. \]

\noin and 

\[ \tilde{E}_{1} \cap q^{-1} (R_{2} )= \left\{ \begin{array}{rl} Y^{4} & =X^{4} -Z^{4} \\

X^{4} U^{4}_{2} & =V^{4}_{2} -W^{4}_{2}  \end{array} \right. \]

Clearly, $\tilde{E}_{1}$ is a smooth complex manifold and 
$q\mid _{\tilde{E}_{1}} :\tilde{E}_{1} \rightarrow \overline{R}$ a holomorphic 
fiber bundle over $\overline{R}$ with fiber $\overline{R}$. Thus 
$\tilde{E}_{1}$ is hyperbolic [1]. The mapping $Q:(X,Y,Z, U_{k}, V_{k} , 
W_{k} )\rightarrow (X,Y,Z, YU_{k} , V_{k} , W_{k} )$ makes 
the pair $(\tilde{E}_{1} , Q\mid _{\tilde{E}_{1}} )$ a normalization of $E_{1}$. \vs

\begin{center} {\bf References} \end{center}

\begin{enumerate}

\item P. Kiernan, {\em Some results in connection with hyperbolic 
manifolds}, Proc. Amer. Math. Soc. {\bf 25} (1970), 588-592. 

\item Sh. Kobayashi, {\em Hyperbolic manifolds and holomorphic mappings}, 
Marcel Dekker, Inc., New York, 1970. 

\item M.H. Kwack, {\em Generalization of big Picard theorem}, Ann. Math. 
{\bf 90} (1969), 9-22. 

\end{enumerate}

\end{document}